# COMPARISON OF TWO MODELS TO PREDICT VERTEBRAL FAILURE LOADS ON THE SAME EXPERIMENTAL DATASET

**V. Allard (1, 2), C. Heidsieck (3), F. Bermond (2), C. Confavreux (1, 4), C. Travert (3,5), L. Gajny, 3, W. Skalli (3), D. Mitton (2), H. Follet (1)**

*1. Univ Lyon, UCBL, INSERM, LYOS 1033, France; 2. Univ Lyon, Univ Eiffel, UCBL, LBMC UMR_T 9406, France; 3. Arts et Métiers Institute of Technologie, Institut de Biomécanique Humaine Georges Charpak, France ; 4. Service de rhumatologie Sud, Hospices Civils de Lyon, France, 5 CHU Pitié Salpêtrière France*

**Abstract.**

Clinical use of finite element analysis requires validation and reproducibility studies. The current study compared two models of vertebral bodies including endplates, on the same experimental dataset and evaluated the influence of the operator on the failure load. Models used were strongly correlated ($R^2=0.91$). The intra-operator reproducibility was 6.4% and 3.5 % for each model. Both simulated results were close to experimental results. The differences in performance could be associated to the differences in segmentation process, mesh (hexahedral vs tetrahedral), material representation and failure criteria. Linear analysis did not decrease model accuracy. Comparison with literature for accuracy and precision shows a wide range of values partly related to the different experimental datasets and the different modelling approaches. Models benchmark using the same experimental dataset are needed to go towards clinical applications.

## 1 Introduction

The models to predict vertebral failure load in literature are numerous (Buckley et al., 2007; Chevalier et al., 2009; Costa et al., 2017; Crawford et al., 2003; Dall'Ara et al., 2022; Giambini et al., 2016; Imai et al., 2006; Johannesdottir et al., 2018; Pahr et al., 2014; Palanca et al., 2020; Prado et al., 2021) but have been rarely applied on the same experimental datasets, restraining comparison between models (Schileo and Taddei, 2021).

Verification and validation (V&V) are processes by which evidence is generated and credibility is thereby established that a computer model yields results with sufficient accuracy for its intended use (ASME Committee (PT60) on Verification and Validation in Computational Solid Mechanics 2006) (Anderson et al., 2007). Since error motivates the need for V&V procedures, it is crucial to understand the types of errors in experimental and computational studies (Anderson et al., 2007).

Clinical use of finite element analysis requires a well-defined process associated to model credibility, following the recommendations of Viceconti et al. (Viceconti et al., 2021).

In this context, the aims of this study are to compare two models of vertebral bodies including endplates, on the same experimental dataset and to evaluate the influence of the operator on the failure load.

## 2 Materials and Methods

### 2.1 Specimens

The experiments were obtained in a previous study (Choisne et al., 2018).
Eleven lumbar spine segments (L1–L3) from cadaveric specimens were considered in this study (5 males and 6 females, age: 82 years± 7 ranging from 61 years to 87 years). Donors were fresh cadavers and no exclusion criteria was specified. Medical history was unveiled but all donors and vertebrae were screened by a trained clinician to ensure there was no evidence of former surgical procedures or spine conditions (e.g. spine curve deviation, suspicion of bone tumor tissue, compression or traumatic fracture). A total of 28 vertebrae were included (8 L1, 11 L2 and 9 L3), after exclusion of vertebrae presenting anomalies.

### 2.2 Imaging

The data acquisition was obtained in a previous study (Choisne et al., 2018). Quantitative Computed Tomography (qCT) scans of the lumbar spines were performed on a Scanner ICT 256 (Philips Healthcare, Cleveland, OH, USA) with the following settings; X-ray tube voltage and current: 120 kV, 1489 mA/s, reconstruction matrix: 512×512, field of view: 250×250mm and a voxel size of 0.39×0.39×0.33 mm. A calibration phantom (QRMESP, QRM GmbH, Germany) was used to map gray scale values to bone mineral density.

### 2.3 Mechanical testing

The mechanical data were obtained in a previous study (Choisne et al., 2018). The 28 vertebrae were extracted from the spinal segment by cutting the intervertebral disc with a scalpel and disarticulating the posterior facets. The vertebrae were then cleaned from all soft tissue and the posterior elements were transected at the pedicles junction to the anterior body. The vertebral endplates were embedded in polymethyl methacrylate (PMMA) for parallelism to ensure uniform loading conditions. Anterior compressive tests were conducted using a quasi-static traction-compression Instron 5500 R device (Instron Ltd., High Wycombe, UK) in combination with a 10 kN load cell and a spherical seating loading platen (Instron Ltd., High Wycombe, UK) allowing for uniform load application onto the upper surface of the specimen (Fig. 1A). The center of rotation of the spherical seating loading platen was aligned with the anterior third of the vertebral body to apply anterior compressive load to the vertebra. The anterior third of the vertebral body was measured using EOS® radiographs of the vertebral body embedded in the PMMA layer before the destructive compression test. The specimen was bolted to the inferior test platen with the mold used for embedding the vertebral body. To avoid shear forces, the inferior test platen was mounted on an antero-posterior sliding connection and an initial compressive preload of 50 N was applied. Specimens were preconditioned by applying 10 cycles of compressive force between 100 N and 250 N at 0.1 Hz followed by destructive compression test at 1 mm/min until ultimate force was achieved (Buckley et al., 2007). Experimental vertebral strength was defined as the ultimate load achieved.

### 2.4 Finite elements models

Two FE models were considered in this study:

**Ensam Model (Choisne et al., 2018)**. A qCT-based finite element model was built from vertebral geometry obtained by a semi-automatic segmentation method (Le Pennec et al., 2014). A hexahedral mesh of the vertebral body (Fig. 1B) was generated from this geometry using a multiblock meshing program wrote in C++ (Grosland et al., 2009). Briefly, the multiblock meshing technique consists in multiple building blocks composed of meshing seeding arranged in rows, columns and layers. The mesh seeds are then projected on the vertebral surface and morphed to each vertebral surface as nodes to lay the foundation for the FE mesh (Grosland et al., 2009). In this 8304-element mesh the average element size was controlled to range between 1mm and 1.5 mm. All meshes were generated with the same topology for each vertebral level allowing the same element to be located approximately at the same position in the vertebra. Convergence analysis was performed to determine the ideal number of elements needed (Choisne et al., 2017). Once the mesh generated, the average BMD of a single finite element was defined on the basis of the qCT scan voxels that fall inside the element. BMD value for each element was converted to linear elastic mechanical properties from an experimental relationship between BMD and elastic modulus (Choisne et al., 2017; Kopperdahl et al., 2002) as shown in Eq. (1). The Poisson ratio, ν, was set to 0.4 (Imai et al., 2006).

$$E(MPa) = 3230 \, BMD \left(\frac{g}{cm^3}\right) - 34.7 \quad (1)$$

*Boundary condition:* Previously described boundary conditions and failure criterion (Sapin-de Brosses et al., 2012) were considered to build the qCT-based model in order to reproduce the experimental testing procedure. Briefly, all models were virtually added layers of PMMA (about 0.5–1 cm thick, E=2500 MPa, ν=0.3) to both vertebral endplates. Lower nodes of the lower PMMA layer were constrained in all degrees of freedom. Anterior compressive load was applied to the upper PMMA layer joined by rigid elements to a node located at the anterior third of the vertebra. Simulations were run on ANSYS software with a linear resolution (ANSYS Inc., Canonsburg, PA, USA). Vertebral strength was defined when a contiguous region of 1cm3 of elements reached 1.5% deformation (Sapin-de Brosses et al., 2012).

**Lyon Model developed by LYOS and LBMC** (Allard et al., 2021). The qCT images were down sampled to 984μm voxel size and then a finite element model was built from vertebral geometry obtained by a semi-automatic method using 3D Slicer (Fedorov 2012, Slicer 4.10.2, https://commonfund.nih.gov/bioinformatics). Briefly, different tools were used to select the vertebral corpse and create a stereolithography (stl) files (threshold, island to delete the noise despeckle, ROI shrink-wrap). From the stl file, a 1 mm$^3$ quadratic tetrahedron mesh (10-nodes) was created (around 50 000 elements for each vertebral corpse). Using the hydroxyapatite phantom composed of several known densities scanned at the same time, specific density from averaged grey levels of each element was assigned to each element using a custom Python script (QCTMA) and converted to Young's modulus using the same relationship in eq (1) from Kopperdahl (Kopperdahl et al., 2002). For simulation purposes, element properties were gathered materials with a fixed step of 10 MPa between Young's modulus of materials (Zannoni et al., 1999). Modulus values resulting from low density regions were set to 100 MPa (Crawford et al., 2003). The Poisson ratio, ν, was set to 0.3.
Nonlinear constitutive law (perfectly elasto-plastic) was considered (ANSYS Inc. v19R1, Canonsburg, PA, USA), with a yield strain of 0.7% (Crawford et al., 2003; Kopperdahl et al., 2002; Pistoia et al., 2002) and a value of compressive strength, defined as the total reaction force generated at an imposed overall vertebral height deformation equivalent to 1.9% strain (Wang et al., 2012).

*Boundary condition:* Before applying boundary conditions, endplates were detected using a custom python script to define the axial direction of the vertebral body.
The displacement and rotation of the surface nodes on the inferior endplate was constrained in all directions, while a displacement load in the z-direction was applied incrementally on the surface nodes on the superior endplate, using a remote point (imposed displacement, rotation allowed) representing a remote ball joint. Remote points are used to simplify the behavior and kinematics of certain portions of the geometry with a single point. They help the user in defining remote boundary conditions like remote displacement and remote force along with elements like joints, springs, beam connections, point masses, moment loads, etc., which simplify the model by reducing the complexity. Anterior compressive load was then applied to the upper endplate joined by rigid elements to a node located at the anterior third of the vertebra, based on the same coordinates as in Ensam Model.

### 2.5 Influence of the operator

According to (Glüer et al., 1995), to minimize the precision error to a 90% confidence interval while evaluating the operator influence, a minimum of 2 observations must be done on at least 27 cases. Hence, to evaluate the inter-operator influence, two operators already fully trained to the model building process applied each model on the entire dataset. To evaluate the intra-operator influence, each analysis was done twice by one operator.
For each model, a protocol was given before applying the finite element analysis, containing specifically:
- phantom calibration equation
- experimental boundary conditions (position, orientation, rate)
To ensure minimum influence, samples were anonymized by the coordinator before they were given to the participants). Coordinator gathered simulated failure load for results analyses.

Results are calculated on the numerical failure load.

$$\text{Intra operator relative difference (\%)} = abs\left(\frac{Trial_2 - Trial_1}{\frac{Trial_2 + Trial_1}{2}} \times 100\right) \quad (2)$$

## 2.6 Statistical analysis

Results are presented as mean ± standard deviation. Accuracy and precision are two measures of observational error. Accuracy is how close a given set of measurements (observations or readings) are to their true value, while precision is how close the measurements are to each other. Here, accuracy is the mean of the difference between simulated and experimental failure load, and precision is the standard deviation of the difference between simulated and experimental failure load.

## 3 Results

Failure loads for experimental data and numerical simulations for both models, all operators and trials are given in Table 1. Experimental failure loads are 3120 ± 1595 N.
Boxplot of experimental data versus numerical results for one operator are shown in figure 2a. ENSAM and Lyon Model results are strongly correlated ($R^2$=0.91) (Fig2b). Bland & Altman graphs (Mean of the difference between simulated and experimental failure load) vs mean between simulated and experimental data is given Figure 3 for both models.
Comparison of the simulated results with the experiments is reported in (Table 2), in terms of accuracy (mean of the difference between simulated and experimental failure load) precision (SD of this difference) and correlation. The mean value of the difference for each trial is given in Table 2.
Each model results are close to the experiments (Table 1).
Results between trials are highly correlated (Table 2). Intra-operator differences in numerical failure load are low (ENSAM: 6.4 ± 6.2%; Lyon: 3.5 ± 2.1 %).

## 4 Discussion and Conclusion

The goals of this study were to compare two models of vertebral bodies including endplates to predict failure load on the same experimental dataset and to evaluate the influence of the operator on the failure load. Both simulated results were strongly correlated (R2=0.91) and were close to experimental results. The differences in performance could result from differences in segmentation process, mesh (hexahedral vs tetrahedral), material representation and failure criteria. Linear analysis did not decrease model accuracy. Comparison with literature for vertebral body with endplates for accuracy and precision is given in figure 4. For a single vertebrae and comparing simulation and experimental results, accuracy range from -1300N (Crawford et al., 2003), -492N (Wang et al., 2012), 68N (Buckley et al., 2007) and -66N (Choisne et al., 2018), and 216 (Ensam Model), 523N (Lyon Model). Precision range from 950N (Crawford et al., 2003), 880N (Wang et al., 2012), 677N (Buckley et al., 2007) and 396N (Choisne et al., 2018), and 340 (Ensam Model), 482N (Lyon Model). Those results show a wide range due in particular to the different experimental datasets.
Intra-operator reproducibility of the simulated failure load could be improved by automation of the segmentation process.
Comparison of models on the same dataset and operator influence are steps needed to assess the credibility of models for future clinical applications.

**Acknowledgement.**
The authors thank the ParisTech BiomecAM chair (COVEA and Société Générale) for their financial support.

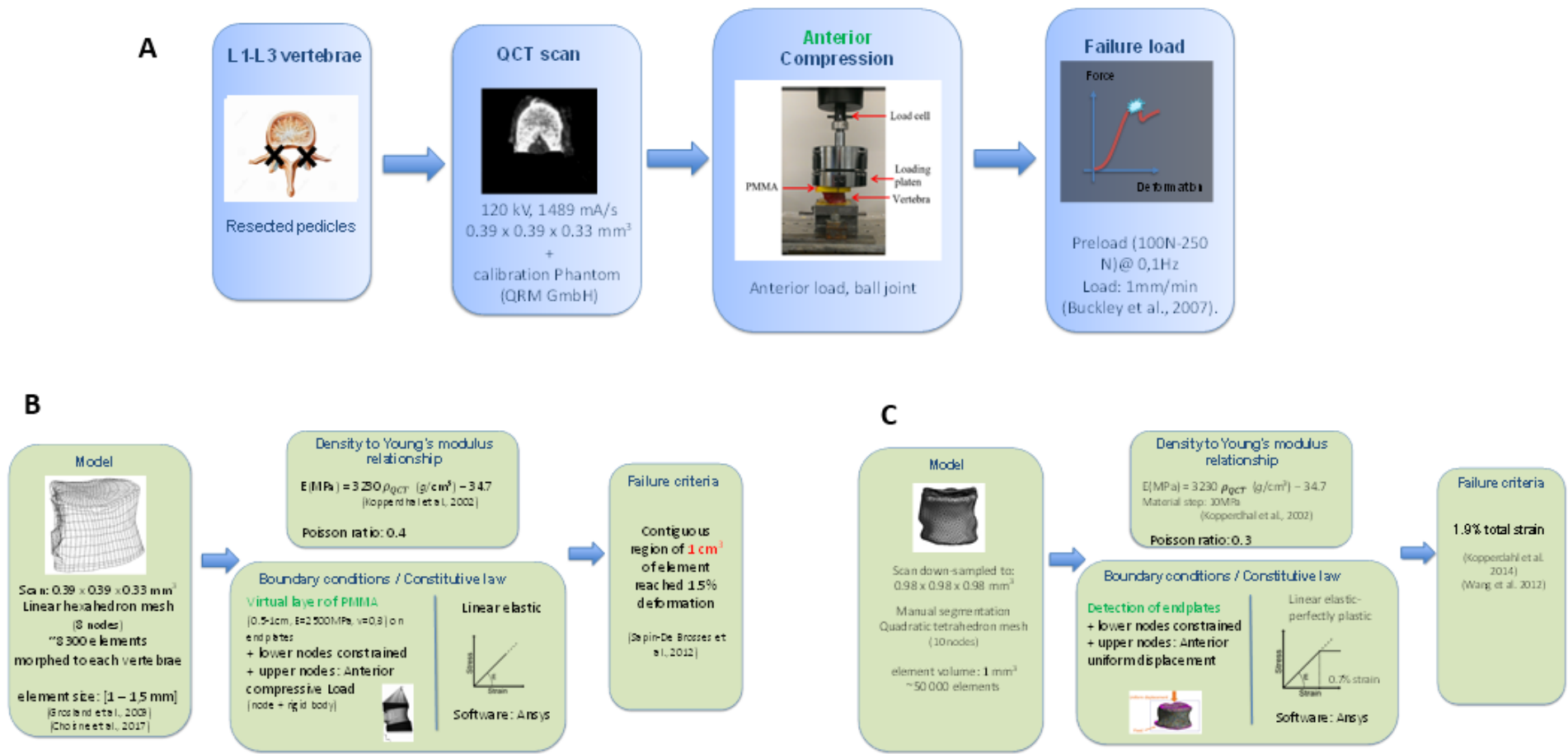

**Figure 1.** A. Experimental Design, adapted from Choisne et al. (Choisne et al. JMBBM 87, 190-196 (2018), B. Model from ENSAM Paris, C: Model from Lyon.

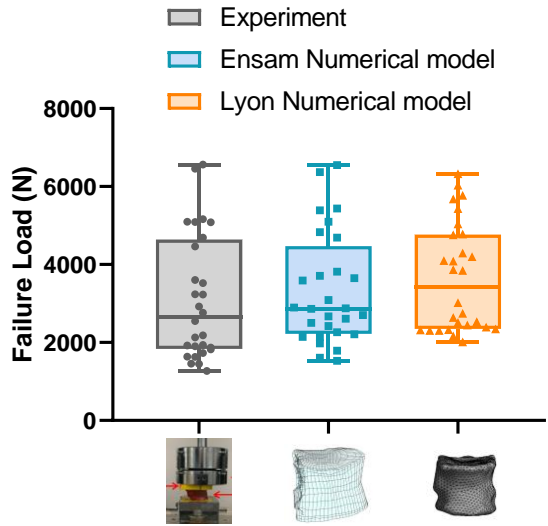 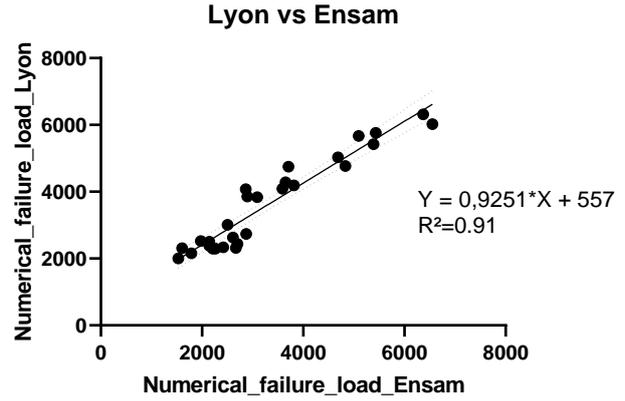

A                                                                B

**Figure 2.** Failure load obtained experimentally and numerically for one operator for each model (for Ensam: operator 1, trial 1, for Lyon: operator 3, trial 1).

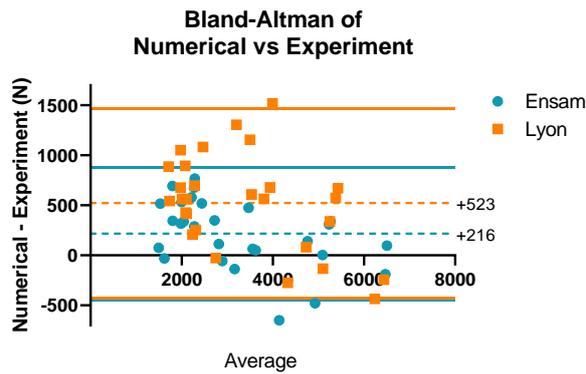

**Figure 3.** Bland- Altman of numerical vs Experiment for one operator for each model (operator 1, trial1 for Ensam, operator 3, trial 1 for Lyon).

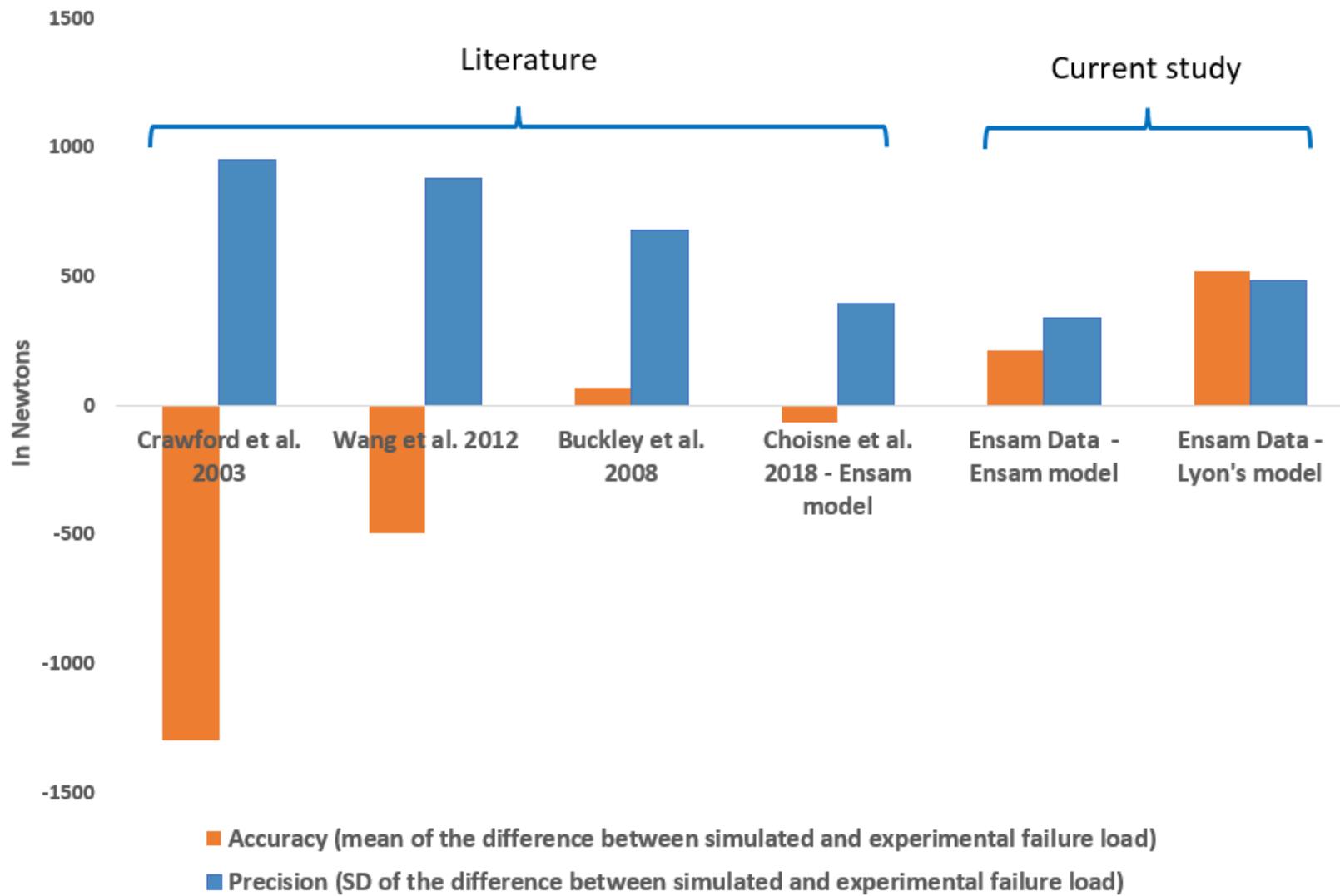

**Figure 4.** Comparison with literature for vertebral body with endplates

**Table 1.** Experimental failure load, and numerical failure load for the different models, operators and tials

| Number | ID | Level | Experimental Failure load (N) | Numerical failure load (N) | | | | |
|---|---|---|---|---|---|---|---|---|
| | | | | Ensam Operator 1 | | Ensam Operator 2 | Lyon Operator 3 | |
| | | | | Trial 1 | Trial2 | Trial 1 | Trial 1 | Trial2 |
| 1 | 438 | L1 | 1935 | 2612 | 2650 | 1954 | 2632 | 2765 |
| 1 | 438 | L2 | 1451 | 2145 | 2334 | 1887 | 2502 | 2641 |
| 1 | 438 | L3 | 1458 | 1533 | 1499 | 1647 | 2001 | 2143 |
| 2 | 421 | L2 | 2133 | 2421 | 2237 | 2129 | 2339 | 2414 |
| 2 | 421 | L3 | 1630 | 1974 | 2139 | 2026 | 2524 | 2621 |
| 3 | 471 | L1 | 3229 | 3092 | 3191 | 3353 | 3839 | 3661 |
| 3 | 471 | L2 | 3233 | 3708 | 3086 | 3523 | 4753 | 4724 |
| 3 | 471 | L3 | 2922 | 2864 | 2726 | 2803 | 4076 | 4034 |
| 4 | 25 | L2 | 6456 | 6552 | 7084 | 6331 | 6021 | 5673 |
| 4 | 25 | L3 | 6561 | 6373 | 6497 | 6094 | 6318 | 6409 |
| 5 | 433 | L1 | 1730 | 2264 | 2140 | 1643 | 2295 | 2376 |
| 5 | 433 | L2 | 1902 | 2670 | 2240 | 2096 | 2317 | 2406 |
| 5 | 433 | L3 | 1825 | 2142 | 2095 | 1819 | 2385 | 2461 |
| 6 | 41 | L1 | 5081 | 5389 | 4684 | 5096 | 5420 | 5679 |
| 6 | 41 | L2 | 5092 | 5094 | 4731 | 4642 | 5667 | 5786 |
| 7 | 383 | L1 | 4465 | 3816 | 3835 | 3668 | 4190 | 4310 |
| 7 | 383 | L2 | 5165 | 4687 | 4698 | 4349 | 5033 | 5125 |
| 7 | 383 | L3 | 4690 | 4832 | 4773 | 4084 | 4771 | 4910 |
| 8 | 484 | L1 | 2182 | 2700 | 2518 | 2499 | 2434 | 2618 |
| 8 | 484 | L2 | 2761 | 2874 | 2862 | 2693 | 2733 | 2944 |
| 8 | 484 | L3 | 1926 | 2505 | 2442 | 2017 | 3009 | 3226 |
| 9 | 448 | L1 | 1270 | 1788 | 1687 | 1929 | 2158 | 2198 |
| 9 | 448 | L2 | 1871 | 2209 | 1758 | 2102 | 2299 | 2298 |
| 9 | 448 | L3 | 1637 | 1607 | 1645 | 1275 | 2313 | 2381 |
| 10 | 436 | L1 | 5095 | 5434 | 4584 | 4246 | 5766 | 5970 |
| 10 | 436 | L2 | 3600 | 3652 | 3662 | 3369 | 4280 | 4440 |
| 10 | 436 | L3 | 3526 | 3590 | 3856 | 3818 | 4090 | 4192 |
| 11 | 493 | L2 | 2548 | 2897 | 2897 | 2782 | 3853 | 3857 |
| | | Mean | 3120 | 3337 | 3234 | 3067 | 3643 | 3724 |
| | | SD | 1595 | 1430 | 1436 | 1362 | 1385 | 1369 |

**Table 2.** Difference Numerical failure load and Experimental failure load, for the different models, operators and trials, and determination coefficient with experimental data

|  | Difference Numerical-Experimental failure load (N) | R² |
|---|---|---|
| Ensam Operator 1- Trial 1 | 216 ± 340 | 0.96 |
| Ensam Operator 1- Trial 2 | 113 ± 385 | 0.94 |
| Ensam Operator 1- Mean Trial 1&2 | 165 ± 331 | 0.96 |
| Ensam Operator 2- Trial 1 | -54 ± 395 | 0.95 |
| Ensam Operators- Mean Trials | 56 ± 337 | 0.97 |
| Ensam intra Operator | 103 ± 298 | 0.96 |
| Ensam inter Operator | 218 ± 279 |  |
|  |  |  |
| Lyon Operator 3- Trial 1 | 523 ± 482 | 0.92 |
| Lyon Operator 3- Trial 2 | 603 ± 504 | 0.91 |
| Lyon Operator 3- Mean Trials | 563 ± 489 | 0.92 |
| Lyon intra Operator | 80 ± 123 | 0.99 |